\documentclass[11pt]{article}

\usepackage{amsmath}
\usepackage{amssymb}
\usepackage{graphicx}

{\catcode `\@=11 \global\let\AddToReset=\@addtoreset}
\AddToReset{equation}{section}

\newtheorem{theorem}{Theorem}[section]
\newtheorem{lemma}{\bf Lemma}[section]

\newtheorem{@definition}{\sc Definition}[section]

\newtheorem{@remark}{\sc Remark}[section]

\newtheorem{@example}{\sc Example}[section]

\def\R{\mathbb{R}}
\def\Z{\mathbb{Z}}

\def\d{\mathrm d}

\def\E{\mathrm E}
\def\P{\Pr}

\def\1{{\bf 1}}
\def\cov{\operatorname{cov}}
\def\var{\operatorname{var}}
\newcommand{\noi}{\noindent}
\newcommand {\nn}{\nonumber}

\newtheorem{thm}{Theorem}[section]
\newtheorem{cor}[thm]{Corollary}

\newtheorem{prop}[thm]{Proposition}

\newtheorem{rem}{Remark}[section]

\begin{document}

\title{Nonparametric estimation 
of the distribution of the autoregressive \\ coefficient
 from panel random-coefficient AR(1) data }
\author{R.\ Leipus$^{1,2}$, \ A.\ Philippe$^{3,4}$, \ V.\ Pilipauskait\.e$^{2,3}$\footnote{Corresponding author. Email address: vytaute.pilipauskaite@gmail.com.}, \  D.\ Surgailis$^2$  }
\date{\today \\ \small
\vskip.2cm
$^1$Vilnius University, Faculty of Mathematics and Informatics, Naugarduko 24, LT-03225 Vilnius, Lithuania\\
$^2$Vilnius University, Institute of Mathematics and Informatics, Akademijos 4, LT-08663 Vilnius, Lithuania\\
$^3$Universit\'{e} de Nantes, Laboratoire de Math\'{e}matiques Jean Leray, 44322 Nantes Cedex 3,
France
\\
$^4$ANJA INRIA Rennes Bretagne Atlantique}

\maketitle

\begin{abstract}
We discuss nonparametric estimation of the distribution function $G(x)$ of the autoregressive coefficient $a \in (-1,1)$ from a panel of $N$ random-coefficient AR(1) data, each of length $n$, by the empirical distribution function of lag~1 sample autocorrelations of individual AR(1) processes.
Consistency and asymptotic normality of the empirical distribution function
and a class of kernel density estimators is established under some regularity conditions on $G(x)$ as $N$ and $n$ increase to infinity. 
The Kolmogorov-Smirnov goodness-of-fit test for simple and composite hypotheses of Beta distributed $a$ is discussed.
A simulation study for goodness-of-fit testing compares the finite-sample performance of our nonparametric estimator to the performance of its parametric analogue discussed in \cite{bib:BSG10}.

\end{abstract}

\medskip


\section{Introduction}

Panel data can describe a large population of heterogeneous units/agents which evolve over time, e.g., households, firms, industries, countries, stock market indices.
In this paper we consider a panel  where  each individual unit
 evolves over time according to order-one random coefficient autoregressive model (RCAR(1)).
 It is well known that aggregation of specific RCAR(1) models  can explain long memory phenomenon, which is often empirically observed in economic time series (see \cite{bib:GR} for instance). More precisely,
consider a panel  $\{ X_i (t), \, t=1,\ldots, n, \, i = 1, \ldots, N\}$,
where each $X_i = \{X_i(t), t \in \Z\}$ is an RCAR(1) process with $(0,\sigma^2)$ noise and random coefficient $a_i\in (-1,1)$, whose autocovariance
\begin{eqnarray}\label{covXi}
\E X_i(0)X_i(t)&=&\sigma^2 \int_{-1}^1 \frac{x^{|t|}}{1- x^2} \,{\d}G(x)   \label{covX}
\end{eqnarray}
is determined by the distribution function $G(x) = \P(a \le x)$ of the autoregressive coefficient.
Granger \cite{bib:GR} showed, for
a specific Beta-type distribution $G(x)$, that the {\it contemporaneous} aggregation of independent processes
$\{ X_i (t) \}$, $i = 1, \ldots, N$,  results in a stationary Gaussian long memory process $\{ {\cal X} (t) \}$, i.e.,
\begin{equation} \label{indaggre}
N^{-1/2} \sum_{i=1}^N X_i (t) \ \to_{\rm fdd} \ {\cal X} (t)
\quad \text{as } N \to \infty,
\end{equation}
where the autocovariance
$\E {\cal X} (0) {\cal X} (t) = \E X_1 (0) X_1 (t)$ decays slowly as $t \to \infty$ so that $\sum_{t \in \Z} | \E {\cal X} (0) {\cal X} (t)  | = \infty$.

A natural statistical problem is recovering
the distribution $G(x)$ (the frequency of $a$ across the population of individual AR(1) `microagents')
from the aggregated sample $\{ {\cal X}(t), \, t=1,\ldots,n\}$.  This problem was treated in \cite{bib:CLP10, bib:CHO06, bib:LOPV06}.
Some related results were obtained in \cite{bib:CE, bib:HL2009, jirak:2013}.  Albeit nonparametric,
the estimators in \cite{bib:CLP10, bib:LOPV06} involve an expansion of the density $g = G'$ in an orthogonal polynomial basis and
are sensitive to the choice of the tuning parameter (the number of polynomials), being limited in practice to very smooth densities $g$.
The last difficulty in estimation of $G$ from aggregated data  is not surprising
due to the fact that aggregation {\it per se} inflicts  a considerable loss of information about the evolution of individual `micro-agents'.

Clearly, if the available data comprises  evolutions  $\{ X_i (t), \, t=1,\ldots,n\}, i=1, \ldots, N$, of all $N$ individual  `micro-agents' (the panel data),
we may expect a much more accurate estimate of $G$.
Robinson \cite{bib:ROB78} constructed an estimator for the moments of $G$ using sample autocovariances of $X_i$
and derived its asymptotic properties  as $N \to \infty$, whereas the length $n$ of each sample remains fixed.
Beran et al. \cite{bib:BSG10} discussed estimation of two-parameter Beta densities $g$ from panel AR(1) data using maximum likelihood
estimators with  unobservable $a_i$ replaced by sample lag 1  autocorrelation coefficient of $X_i (1), \ldots, X_i (n)$ (see Section \ref{s:sim}), and
derived the asymptotic normality together with some other properties of the estimators as $N$ and $n$ tend to infinity.

The present paper studies nonparametric estimation of $G$ from panel random-coefficient AR(1) data using the
empirical distribution function:
\begin{equation} \label{hatG0}
\widehat G_{N,n}(x)\ := \ \frac{1}{N} \sum_{i=1}^N \1( \widehat a_{i,n} \le x),
\quad x \in \R,
\end{equation}
where $\widehat a_{i,n}$ is the lag~1 sample autocorrelation coefficient of $X_i$, $i=1,\ldots, N$
(see \eqref{hatan} below). We also discuss kernel estimation of the density $g(x) = G'(x)$ based on smoothed version of \eqref{hatG0}.
We assume that individual AR(1) processes $X_i$ are driven by identically distributed shocks containing both common and idiosyncratic (independent) components.
Consistency and asymptotic normality as $N, n \to \infty $
of the above estimators are derived under some regularity conditions on $G(x)$. Our results can be applied to test goodness-of-fit
of the
distribution $G(x)$ to a given hypothesized distribution (e.g., a Beta distribution) using the Kolmogorov-Smirnov statistic, and
to construct  confidence intervals for $G(x)$ or $g(x)$.

The paper is organized as follows. Section 2 obtains the  rate of convergence of the sample  autocorrelation coefficient $\widehat a_{i,n}$ to $a_i$,
in probability,
the result of independent interest. Section 3 discusses the weak convergence of the empirical process in \eqref{hatG0}
to a generalized Brownian bridge. The Kolmogorov-Smirnov goodness-of-fit test
for simple and composite hypotheses of Beta distributed $a$ is discussed in Section 4.
In Section 5 we study kernel density estimators of $g(x)$. We
show that these estimates are asymptotically normally distributed and
their mean integrated square error tends to zero. A simulation study of Section 6 compares the empirical performance
of  \eqref{hatG0} and the parametric estimator of \cite{bib:BSG10} to the goodness-of-fit testing for $G(x)$ under  null
Beta distribution. The proofs of auxiliary statements
can be found in the Appendix.

In what follows, $C$ stands for a positive constant whose precise value is unimportant and which may change from line to line. We write  $\to_{\rm p}, \, \to_{\rm d}, \, \to_{\rm fdd}$ for the convergence in probability and the convergence of (finite-dimensional) distributions respectively, whereas $\Rightarrow$ denotes the weak convergence in the space $D[-1,1]$ with the supremum metric.

\section{Estimation of random autoregressive coefficient}

Consider an RCAR(1) process
\begin{equation}\label{def:AR1}
X(t) \ = \ a X(t-1) + \zeta(t), \quad t \in \Z,
\end{equation}
where innovations $\{ \zeta (t) \}$ admit the following decomposition:
\begin{equation}\label{def:fs}
\zeta (t) \ = \ b \eta (t) + c \xi (t), \quad t \in \Z,
\end{equation}
where random sequences $\{\eta (t)\}$, $\{ \xi(t) \}$ and random coefficients
$a$, $b$, $c$  satisfy the following conditions:
\smallskip

\noi {\bf Assumption A$_1$}\ \ $\{ \eta (t) \}$ are independent identically distributed (i.i.d.) random variables (r.v.s) with $\E \eta (0) = 0$, $\E \eta^2 (0) = 1$, $\E | \eta (0)|^{2p} < \infty$ for some $p> 1$.
\smallskip

\noi {\bf Assumption A$_2$}\ \ $\{ \xi (t) \}$ are i.i.d.\ r.v.s with $\E \xi (0) = 0$, $\E \xi^2 (0) = 1$, $\E | \xi (0) |^{2p} < \infty$ for the same $p$ as
in A$_1$.
\smallskip

\noi {\bf Assumption A$_3$}\ \ $b$ and $c$ are possibly dependent r.v.s such that  $\P ( b^2 + c^2 > 0 ) = 1$ and $\E b^2 < \infty$, $\E c^2 < \infty$.\smallskip

\noi {\bf Assumption A$_4$}\ \ $a\in (-1,1)$ is a r.v.\ with a distribution function (d.f.)
$G(x) := \P(a\le x)$ supported on $[-1,1]$ and satisfying
\begin{eqnarray}\label{cond:Ea}
\E \Big( \frac{1}{1-|a|} \Big) &=& \int_{-1}^1 \frac{{\rm d} G(x)}{1-|x|} \ <\ \infty.
\end{eqnarray}

\noi {\bf Assumption A$_5$}\ \ $a$, $\{ \eta (t) \}$, $\{ \xi (t) \}$ and the vector $(b,c)^\top$ are mutually independent.
\smallskip

\begin{rem} {\rm In the context of panel observations (see \eqref{def:multple_AR1} below), $\{\eta(t)\}$ is the common component
and $\{ \xi (t) \}$ is the idiosyncratic component of shocks. The innovation process $\{\zeta (t) \}$ in \eqref{def:fs} is i.i.d.\ if
the coefficients $b $ and $c$ are nonrandom. In the general case $\{\zeta (t) \}$ is a dependent and uncorrelated stationary process
with $\E \zeta (0) = 0$,
$\E \zeta^2 (0) = \E b^2 + \E c^2$, $\E \zeta (0) \zeta (t) = 0$, $t \not = 0$. }

\end{rem}

Under conditions A$_1$--A$_5$, a unique strictly stationary solution of \eqref{def:AR1} with finite variance exists and is written as
\begin{eqnarray}\label{Xma}
 X(t) \ = \ \sum_{s\le t} a^{t-s} \zeta(s), \quad t \in \Z.
\end{eqnarray}
Clearly, $\E X(t) = 0$ and $\E X^2(t) = \E \zeta^2 (0) \E (1-a^2)^{-1} < \infty$.
Note that \eqref{cond:Ea} is equivalent to
$$
\E \Big( \frac{1}{1-|a|^p} \Big) \ < \ \infty, \quad 1 < p \le 2,
$$
since $1-|a| \le 1 - |a|^p \le  2(1 - |a|)$ for $a \in (-1,1)$.

For an observed sample $X(1),\ldots,X(n)$ from the stationary process
in \eqref{Xma},
define the sample mean
$\bar X_n := n^{-1} \sum_{t=1}^n X(t)$ and
the sample lag~1 autocorrelation coefficient
\begin{eqnarray}\label{han}
\widehat{a}_n \ := \ \frac{\sum_{t=1}^{n-1} (X(t) - \bar X_n)  (X(t+1)- \bar X_n)}{\sum_{t=1}^{n} (X(t)- \bar X_n)^2}.
\end{eqnarray}
Note the estimator $\widehat{a}_n$ in \eqref{han}  does not exceed 1 a.s.\ in absolute value
by the Cauchy-Schwarz inequality. Moreover, it
is invariant to shift and scale transformations of $\{ X(t) \}$ in \eqref{def:AR1}, i.e., we can replace $\{  X(t)\}$ by $\{ \rho X(t) + \mu \}$ with some (unknown) $\mu \in \R$ and $\rho > 0$. 

\begin{prop}\label{prop:1} Under Assumptions A$_1$--A$_5$, for any $0 < \gamma < 1$ and $n \ge 1$, it holds
\begin{equation*}
\P ( |\widehat{a}_n - a| > \gamma) \ \le \ C (n^{-((p/2)\wedge (p-1)) } \gamma^{-p} + n^{-1} ),
\end{equation*}
with $C > 0$ independent of $n, \gamma$.
\end{prop}

\noi{\sc Proof.} See Appendix.
\medskip

Assume now that the d.f.\ $G(x)=\P(a\le x)$ satisfies the following H\"older condition:
\smallskip

\noi \textbf{Assumption A$_6$}\ \ There exist constants $L_G >0$ and $\varrho \in (0,1]$ such that
\begin{equation} \label{cond:GHoe}
 |G(x)-G(y)| \ \le \ L_G |x-y|^\varrho, \quad x,y \in [-1,1].
\end{equation}

Consider the d.f.\  of $\widehat{a}_n$:
\begin{equation}\label{def:cdf_hatan}
G_n (x) \ := \ \P (\widehat{a}_n \le x),\quad x \in \R.
\end{equation}


\begin{cor}\label{cor:1} Let Assumptions A$_1$--A$_6$ hold. Then, as $n \to \infty$,
 $$
 \sup_{x \in [-1,1]}|G_n (x) - G(x)| = O(n^{-\frac{\varrho}{\varrho+p} (\frac{p}{2} \wedge (p-1))}).
 $$
\end{cor}

\noi{\sc Proof.} Denote $\delta_n := \widehat{a}_n - a$. For any (nonrandom) $\gamma>0$ from
\eqref{cond:GHoe} we have
\begin{eqnarray*}
\sup_{x \in [-1,1]}|G_n (x) - G(x)| & = & \sup_{x \in [-1,1]} |\P (a+ \delta_n \le x)- \P (a\le x)|\
\le \  L_G\gamma^\varrho +\P(| \delta_n |>\gamma),
\end{eqnarray*}
implying
\begin{equation*}
\sup_{x \in [-1,1]}| G_n (x) - G(x)| \ \le \ L_G \gamma^\varrho + C(n^{-1} +  n^{-((p/2)\wedge (p-1))} \gamma^{-p})
\end{equation*}
with $C > 0$ independent of $n, \gamma$.
Then the corollary follows from  Proposition \ref{prop:1} by taking $\gamma = \gamma_n  = o(1) $ such
that $\gamma^{\varrho}_n \sim n^{-((p/2)\wedge (p-1))} \gamma^{-p}_n$ and noting that the exponent
$\frac{\varrho}{\varrho+p} (\frac p 2 \wedge (p-1)) < 1$.
 \hfill $\Box$

\section{Asymptotics of the empirical distribution function}

Consider random-coefficient AR(1) processes $\{ X_i(t)\}$, $i=1, 2, \dots$, which are stationary solutions to
\begin{equation}\label{def:multple_AR1}
X_i(t) \ = \ a_i X_i (t-1) + \zeta_i (t), \quad t \in \Z,
\end{equation}
with innovations $\{ \zeta_i (t) \}$ having the same structure as in \eqref{def:fs}:
\begin{equation}\label{def:innov}
\zeta_i (t) \ = \ b_i \eta(t) + c_i \xi_i(t), \quad t \in \Z.
\end{equation}
More precisely, we make the following assumption:

\smallskip

\noi {\bf Assumption B}\ \
$\{ \eta (t) \}$ satisfies A$_1$; $\{ \xi_i (t) \}$, $(b_i, c_i)^\top$, $a_i$, $i=1, 2, \ldots,$ are independent copies of
$\{ \xi (t) \}$, $(b, c)^\top$, $a$, respectively, which satisfy Assumptions A$_2$--A$_6$. (Note that we assume A$_5$ for any $i=1, 2, \ldots$.)

\smallskip

\begin{rem}
{\rm The individual processes $X_i$ have covariance long memory 
if conditions \eqref{cond:Ea} and 
$\int_{-1}^1 |1 - x|^{-2} \d G(x) = \infty$ hold, which is compatible with Assumption~B.  The same is true 
about the limit aggregated process in \eqref{indaggre}  arising when the common component is absent. 
On the other hand, in the presence of the common component, long memory in the limit aggregated process arises when 
the individual processes have infinite variance 
and condition \eqref{cond:Ea} fails, see \cite{bib:PUP09}. 
}
\end{rem}

Define the sample mean $\bar X_{i,n} := n^{-1} \sum_{t=1}^n X_i(t)$, the corresponding sample lag~1 autocorrelation coefficient
\begin{eqnarray} \label{hatan}
\widehat a_{i,n} \ := \ \frac{\sum_{t=1}^{n-1} (X_i(t)-\bar X_{i,n})
(X_i(t+1)-\bar X_{i,n})}{\sum_{t=1}^{n} (X_i(t) - \bar X_{i,n})^2}, \qquad 1 \le i \le N
\end{eqnarray}
and the empirical d.f.
\begin{eqnarray} \label{hatG}
\widehat G_{N,n}(x)\ := \ \frac{1}{N} \sum_{i=1}^N \1( \widehat a_{i,n} \le x),
\quad x \in \R.
\end{eqnarray}
Recall that \eqref{hatG} is  a nonparametric estimate of the d.f.\ $G(x) = \P(a_i \le x)$ from observed panel data $\{X_i(t), \, t=1,\ldots,n, \, i=1,\ldots,N\}$.
In the following theorem we show that $\widehat G_{N,n}(x)$ is an asymptotically unbiased estimator
of $G(x)$, as $n$ and $N$ both tend to infinity, and prove the weak convergence of the corresponding empirical process.

\begin{theorem}\label{thm:1} Assume the panel data model in \eqref{def:multple_AR1}--\eqref{def:innov}.
Let Assumption B hold and $N, n \to \infty$.  Then
\begin{equation}\label{lim:sup_bias}
 \sup_{x\in [-1,1]} |\E\widehat G_{N,n}(x)-G(x)| \ = \  O(n^{- \frac{\varrho}{\varrho + p} ( \frac{p}{2} \wedge (p-1))}).
\end{equation}
If, in addition,
\begin{equation} \label{rateNn}
N = o( n^{\frac{2 \varrho}{\varrho+p} ( \frac{p}{2} \wedge (p-1) )}),
\end{equation}
then
\begin{equation}\label{lim:emp_proc}
N^{1/2} (\widehat G_{N,n}(x) - G(x)) \ \Rightarrow \ W(x), \quad x \in [-1,1],
\end{equation}
where $\{W(x), \, x \in [-1,1]\}$ is a continuous Gaussian process with zero mean and 
$\cov (W(x), W(y)) =  G(x\wedge y) - G(x) G(y)$. 
\end{theorem}

\noi{\sc Proof.} Note $\widehat{a}_{i,n}$, $i = 1, \ldots, N$, are identically distributed, in particular, $\E \widehat{G}_{N,n} (x) = G_n (x)$ with $G_n (x)$ defined in \eqref{def:cdf_hatan}. Hence, \eqref{lim:sup_bias} follows immediately from Corollary \ref{cor:1}.
\smallskip

To prove the second statement of the theorem, we approximate $\widehat G_{N,n}(x)$ by the empirical d.f.
$$
\widehat G_N (x) \ := \ \frac{1}{N} \sum_{i=1}^N \1(a_i \le x), \qquad x \in  [-1,1]
$$
of i.i.d.\ r.v.s $a_i, i=1,\ldots, N$.
We have $N^{1/2} (\widehat G_{N,n}(x) - G(x)) = N^{1/2} (\widehat G_N(x) - G(x)) + D_{N,n} (x)$
with $D_{N,n}(x) := N^{1/2} (\widehat G_{N,n}(x) - \widehat G_N(x))$. Since A$_6$ guarantees the continuity of $G$, it holds
\begin{equation*}
N^{1/2} (\widehat G_N (x) - G(x)) \ \Rightarrow \ W(x), \quad x \in [-1,1]
\end{equation*}
by the classical Donsker theorem. Then \eqref{lim:emp_proc} follows once we prove
$\sup_{x\in [-1,1]}  | D_{N,n} (x) | \to_{\rm p} 0$. By definition,
$$
D_{N,n} (x) \ = \ N^{-1/2} \sum_{i=1}^N
(\1 (a_i + \delta_{i,n} \le x) - \1(a_i \le x)) \ = \ D'_{N,n}(x) - D''_{N,n}(x),$$
where $\delta_{i,n}:= \widehat{a}_{i,n} - a_i$, $i=1,\ldots,N$, and
\begin{eqnarray*}
D'_{N,n}(x) &:=& N^{-1/2} \sum_{i=1}^N \1(x< a_i \le x- \delta_{i,n},\, \delta_{i,n} \le 0), \\
D''_{N,n}(x) &:=& N^{-1/2} \sum_{i=1}^N \1(x - \delta_{i,n} < a_i \le x, \, \delta_{i,n} > 0).
\end{eqnarray*}
For $\gamma > 0$ we have
$$
D'_{N,n}(x) \ \le \ N^{-1/2} \sum_{i=1}^N \1(x< a_i \le x + \gamma) + N^{-1/2} \sum_{i=1}^N \1 ( | \delta_{i,n}|>\gamma) \ =: \
V'_{N}(x) + V''_{N,n}.
$$
(Note that $V''_{N,n}$ does not depend on $x$.) By Proposition \ref{prop:1},
we obtain
$$
\E V''_{N,n}
\ = \ N^{-1/2} \sum_{i=1}^N \P(| \delta_{i,n}| >\gamma)
\ \le \ C N^{1/2} ( n^{- ((p/2) \wedge (p-1))} \gamma^{-p} + n^{-1}),
$$
which tends to $0$
when $\gamma $ is chosen as $\gamma^{\varrho+p} = n^{-((p/2) \wedge (p-1))} \to  0$. Next,
\begin{eqnarray*}
V'_{N}(x) &=& N^{1/2}(\widehat G_N(x+\gamma) -  \widehat G_N(x)) \ = \ N^{1/2}(G(x+\gamma) - G(x)) + U_N (x, x + \gamma],\\
U_N (x, x + \gamma]  &:=& N^{1/2} ( \widehat G_N (x+\gamma) - G (x +\gamma) ) - N^{1/2} ( \widehat G_N (x) - G (x) ).
\end{eqnarray*}
The above choice of $\gamma^{\varrho+p} = n^{-((p/2) \wedge (p-1))}$ implies
$\sup_{x \in [-1,1]} N^{1/2}|G(x+\gamma)- G(x)| = O(N^{1/2} \gamma^\varrho) = o(1)$,
whereas $U_N (x, x + \gamma]$ vanishes in the uniform metric in probability (see Lemma \ref{lemma:2} in Appendix).
 Since $D''_{N,n}(x)$ is analogous to $D'_{N,n}(x)$, this proves the theorem.
 \hfill $\Box$

\begin{rem} {\rm \eqref{rateNn} implies that $n \gg N^{(\varrho + p)/\varrho p} $ asymptotically for $p \ge 2$. Note that $(\varrho + p)/\varrho p > 1 $ and $\lim_{p \to \infty} (\varrho + p)/\varrho p =  1/\varrho $ for any $\varrho \in (0,1]$.
We may conclude that Theorem \ref{thm:1} as well as other results of this  paper apply
to  {\it long} panels with $n $ increasing much faster than $N$, except maybe for the limiting case
$p= \infty $ for $\varrho =1 $. The main reason for this conclusion is that
$a_i$ need to be accurately estimated by \eqref{hatan} in order that 
$\widehat G_{N,n}(x)$ behaves similarly to the 
empirical d.f.\ $\widehat G_{N}(x)$ based on unobserved
autocorrelation coefficients $a_i, 1 \le i \le N$. }

\end{rem}

\section{Goodness-of-fit testing}

Theorem \ref{thm:1} can be used for testing goodness-of-fit. In the case of {\it simple} hypothesis, we test the null
$H_0\colon G=G_0$ vs.\ $H_1\colon G \neq G_0$ with $G_0$ being a certain hypothetical distribution satisfying the H\"older
condition in \eqref{cond:GHoe}. Accordingly,
the corresponding Kolmogorov-Smirnov (KS) test rejecting $H_0$ whenever
\begin{equation}\label{gof}
N^{1/2} \sup_{x \in [-1,1]}  | \widehat{G}_{N,n}(x) - G_0 (x) | > c (\omega)
\end{equation}
has asymptotic size $\omega \in (0,1)$ provided $N, n, G_0$ satisfy the assumptions for \eqref{lim:emp_proc}
in Theorem \ref{thm:1}. (Here, $c(\omega)$ is the upper $\omega$-quantile of the Kolmogorov distribution.)
However, the goodness-of-fit test in \eqref{gof}
requires the knowledge of parameters of the model considered, which is not typically a very realistic situation.
Below, we consider testing  {\it composite} hypothesis
using the Kolmogov-Smirnov statistic with estimated parameters. The parameters will be estimated by
the method of moments.

\smallskip

Write $\mu = (\mu^{(1)}, \ldots, \mu^{(m)})^\top$ and $\widehat{\mu}_{N,n} = (\widehat{\mu}^{(1)}_{N,n}, \ldots, \widehat{\mu}^{(m)}_{N,n} )^\top$, where
$$
\mu^{(u)} := \E a^u = \int_{-1}^{1} x^u \d G (x), \qquad \widehat{\mu}_{N,n}^{(u)} := \frac{1}{N} \sum_{i=1}^{N} (\widehat{a}_{i,n} )^u,
\quad 1 \le u \le m.
$$

\begin{prop}\label{prop:0} Let the panel data model in \eqref{def:multple_AR1}--\eqref{def:innov}
satisfy Assumption B with exception of Assumption~A$_6$. If $N = o(n^{\frac{2}{1+p} (\frac{p}{2} \wedge (p-1))})$ as $N, n \to \infty$, then
\begin{eqnarray} \label{cltM}
N^{1/2} ( \widehat{\mu}_{N,n}  - \mu) \ \to_{\rm d} \ {\mathcal N} (0, \Sigma), \quad \text{where} \quad
\Sigma := \big(\cov (a^u, a^v)\big)_{1 \le u,v \le m}.
\end{eqnarray}
\end{prop}

\noi {\it Proof.} Write
\begin{eqnarray*}
N^{1/2} ( \widehat{\mu}_{N,n}  - \mu) &=&N^{1/2} ( \widehat{\mu}_{N,n}  - \widehat \mu_N) + N^{1/2} ( \widehat{\mu}_N  - \mu),
\end{eqnarray*}
 where $\widehat{\mu}_N := \frac{1}{N} \sum_{i=1}^N (a_i, \ldots, a_i^m)^\top$. We have $N^{1/2} ( \widehat{\mu}_N - \mu ) \to_{\rm d} {\mathcal N} (0, \Sigma)$ as $N \to \infty$ by the multivariate central limit theorem. On the other hand, $N^{1/2} (\widehat{\mu}_{N,n} - \widehat{\mu}_N) \to_{\rm p} 0$
follows from
$\E | \widehat{a}_n^u - a^u | \le C \E | \widehat{a}_n - a | \le C (\gamma + \P ( |\widehat{a}_n - a | > \gamma ))$ and
Proposition~\ref{prop:1} with $\gamma^{1+p} = n^{-((p-1) \wedge (p/2))}$,
proving the proposition. \hfill $\Box$

\begin{rem} {\rm 
Robinson \cite[Theorem 7]{bib:ROB78} discussed a different estimate of
$\mu$, which was proved to be asymptotically normal for fixed $n$ as $N \to \infty$ in contrast to ours.
However, his result holds
in the case of idiosyncratic innovations only and
under stronger assumptions on $G$ than in Proposition \ref{prop:0},
which do not allow for long memory.  }
\end{rem}

Consider testing the composite null hypothesis that $G$ belongs to the family
${\cal G} = \{ G_\theta, \, \theta=(\alpha,\beta)^\top \in (1,\infty)^2 \}$ of
Beta d.f.s versus an alternative $G \not\in {\cal G}$, where
\begin{equation}\label{def:beta}
G_\theta (x) = \frac{1}{B(\alpha, \beta)} \int_0^x t^{\alpha - 1} (1 - t)^{\beta - 1} \d t, \qquad x \in [0,1]
\end{equation}
and ${\rm B}(\alpha, \beta) = \Gamma(\alpha) \Gamma(\beta) / \Gamma(\alpha + \beta)$ is Beta function. 
The $u$th moment of  $G_\theta $ is given by
\begin{equation*}\label{def:mu}
\mu^{(u)} = \int_{0}^{1} x^u \d G_\theta (x) = \prod_{r=0}^{u-1} \frac{\alpha + r}{\alpha + \beta + r}.
\end{equation*}
Parameters $\alpha, \beta $ can be found from the first two moments $\mu = (\mu^{(1)}, \mu^{(2)})^\top$ as
\begin{equation}\label{def:theta}
\alpha = \frac{\mu^{(1)} (\mu^{(1)} - \mu^{(2)})}{\mu^{(2)} - (\mu^{(1)})^2}, \qquad \beta =  \frac{(1 - \mu^{(1)}) (\mu^{(1)} - \mu^{(2)})}{\mu^{(2)} - (\mu^{(1)})^2} .
\end{equation}
The moment-based estimator $\widehat{\theta}_{N,n} := (\widehat{\alpha}_{N,n}, \widehat{\beta}_{N,n})^\top$ of $\theta = (\alpha, \beta)^\top$ is
obtained by replacing $\mu $
in \eqref{def:theta} by its estimator $\widehat{\mu}_{N,n}$.
The consistency and asymptotic normality of this estimator follows by  the Delta method from Proposition~\ref{prop:0}, see Corollary \ref{cor:2} below,
where we need condition $\alpha, \beta > 1$ to satisfy Assumptions A$_4$ and A$_6$.

\begin{cor}\label{cor:2}
Let the panel data model in \eqref{def:multple_AR1}--\eqref{def:innov}
satisfy Assumption B.
Let $G = G_\theta$, $\theta = (\alpha, \beta)^\top$ be a Beta d.f.\ in \eqref{def:beta}, where $\alpha > 1$, $\beta > 1$.
Let $N, n $ increase as in \eqref{rateNn} where $\varrho = 1$.
Then
\begin{equation}\label{clt}
N^{1/2} ( \widehat{\theta}_{N,n} - \theta ) \ \to_{\rm d} \ {\mathcal N} (0, \Lambda_\theta), \qquad \Lambda_\theta := \Delta^{-1} \Sigma (\Delta^{-1})',
\end{equation}
where $\Sigma $ is the $2\times 2$ matrix in \eqref{cltM} and
\begin{eqnarray*}
\Delta  := \partial \mu / \partial \theta   =  \begin{pmatrix}
\partial \mu^{(1)} / \partial \alpha & \partial \mu^{(1)} / \partial \beta \\
\partial \mu^{(2)} / \partial \alpha & \partial \mu^{(2)} / \partial \beta
\end{pmatrix}.
\end{eqnarray*}
Moreover, $\widehat{\theta}_{N,n}$ is asymptotically linear:
\begin{eqnarray}\label{def:l}
N^{1/2} ( \widehat{\theta}_{N,n} - \theta ) = N^{-1/2} \sum_{i=1}^N l_\theta (a_i) + o_{\rm p} (1), \qquad
l_\theta (x) := \Delta^{-1}
(x - \mu^{(1)},
x^2 - \mu^{(2)})^\top,
\end{eqnarray}
where  $\E l_\theta (a) = \int_0^1 l_\theta (x) \d G_\theta (x) = 0$ and $\E l_\theta (a) l_\theta (a)^\top = \int_0^1 l_\theta (x) l_\theta (x)^\top \d G_\theta (x) = \Lambda_\theta$.
\end{cor}

\begin{cor}\label{cor:3}
Let assumptions of Corollary~\ref{cor:2} hold. Then
$$
N^{1/2} (\widehat{G}_{N,n} (x) - G_{\widehat{\theta}_{N,n}} (x) ) \ \Rightarrow \ V_\theta (x), \qquad x \in [0,1],
$$
where $\{ V_\theta (x), \, x \in [0,1] \}$ is a continuous Gaussian process with zero mean and covariance
\begin{eqnarray*}
\cov ( V_\theta (x), V_\theta (y) ) &=& G_\theta (x \wedge y) - G_\theta (x) G_\theta (y)\\
&&- \int_0^x l_\theta (u)^\top \d G_\theta (u) \partial_\theta G_\theta (y) - \int_0^y l_\theta (u)^\top \d G_\theta (u) \partial_\theta G_\theta (x) +
\partial_\theta G_\theta (x)^\top \Lambda_\theta \partial_\theta G_\theta (y),
\end{eqnarray*}
where $\partial_\theta G_\theta (x) := \partial G_\theta (x) / \partial \theta =
\big(\partial G_\theta (x) / \partial \alpha, \partial G_\theta (x) / \partial \beta \big)^\top$, $x \in [0,1]$ and $\Lambda_\theta$ is defined
in \eqref{clt}.
\end{cor}

\noi {\it Proof.} The d.f.\ $G_\theta$ with $\alpha > 1, \beta > 1$ satisfies Assumptions A$_4$ and A$_6$ with $\varrho = 1$.
Recall $\widehat{G}_N (x) := N^{-1} \sum_{i=1}^{N} \1 (a_i \le x)$, $x \in [0,1]$.
Since condition \eqref{rateNn} is satisfied, so $N^{1/2}\sup_{x\in [0,1]} |\widehat{G}_{N,n} (x) - \widehat{G}_N (x) |$ vanishes in probability by Theorem~\ref{thm:1}, whereas the convergence $N^{1/2} (\widehat{G}_N (x) - G_{\widehat{\theta}_{N,n}} (x) ) \Rightarrow V_\theta (x)$, $x \in [0,1]$ follows from \eqref{def:l} using the fact that $ \partial_\theta G_\theta (x), \, x \in [0,1]$ is continuous in $\theta$, see \cite{bib:durb1973} or \cite[Theorem 19.23]{bib:vaart2000}.
\hfill $\Box$
\medskip

With Corollary~\ref{cor:3} in mind, the Kolmogorov-Smirnov test for the composite hypothesis
$G \in {\cal G}$ can be defined as
\begin{eqnarray}\label{KStest}
\sup_{x\in [0,1]} N^{1/2} |\widehat{G}_{N,n} (x) - G_{\widehat{\theta}_{N,n}} (x)|&>&c(\omega,  \widehat{\theta}_{N,n}),
\end{eqnarray}
where $c(\omega,\theta) $ is the upper $\omega$-quantile of the distribution of
$\sup_{x\in [0,1]}  |V_\theta (x)|$:
\begin{equation*}
\P \Big( \sup_{x\in [0,1]}  |V_\theta (x)| > c(\omega,\theta) \Big) = \omega, \qquad \omega \in (0,1).
\end{equation*}
The test in \eqref{KStest}  has correct asymptotic size for any $\omega \in (0,1)$, which follows
from Corollary \ref{cor:3} and
the continuity of the quantile function $c(\omega, \theta)$ in $\theta$, see
\cite[p.~69]{bib:szuc2008}, \cite{bib:vaart2000}.
By writing
$N^{1/2} (\widehat{G}_{N,n} (x) - G_{\widehat{\theta}_{N,n}} (x)) =
N^{1/2} (\widehat{G}_{N,n} (x) - G(x)) + N^{1/2} (G(x) - G_{\widehat{\theta}_{N,n}} (x))$,
it follows that the Kolmogorov-Smirnov statistic on the l.h.s.\ of \eqref{KStest} tends to infinity (in probability) under any
fixed alternative $G \not\in {\cal G}$ which cannot be approximated by a Beta d.f. $G_\theta$ in the uniform
metric, i.e.,
such that $\inf_{\theta}
\sup_{x \in [0,1]} |G(x) - G_\theta (x)| > 0$. Moreover, even under the alternative, we preserve the consistency of $\widehat \mu_{N,n}$, hence $c(\omega, \widehat{\theta}_{N,n})$ being a continuous function of sample moments,
converges in probability to some finite limit. Therefore the test \eqref{KStest} is consistent.

In practice, the evaluation of $c(\omega, \theta)$ requires Monte Carlo approximation which is time-consuming. Alternatively, 
\cite{bib:stut1993,bib:szuc2008} discussed 
parametric bootstrap procedures to produce asymptotically correct critical values. We note that  
the assumptions of \cite[Theorem~1]{bib:szuc2008} are valid for the family of Beta d.f.s and the moment-based estimator of $\theta$
in Corollary~\ref{cor:3}. 
The consistency of the test when using bootstrap critical values follows by a similar argument as in \eqref{KStest}.

\section{Kernel density estimation
}

In this section we assume $G$ has a bounded probability density function $g(x) = G'(x), x \in [-1,1]$, implying
Assumption  A$_6$
with H\"older exponent $\varrho =1$ in \eqref{cond:GHoe}.
It is of interest to estimate $g(x)$ in a nonparametric way from $\widehat{a}_{1,n}, \ldots, \widehat{a}_{N,n}$  \eqref{hatan}.

Consider the kernel density estimator
\begin{eqnarray}\label{def:g_nN}
\widehat{g}_{N,n} (x) \ := \ \frac{1}{N h} \sum_{i=1}^{N} K \Big( \frac{x - \widehat{a}_{i,n}} {h} \Big),
\quad x \in \R,
\end{eqnarray}
where $K$ is a kernel, satisfying Assumption A$_7$ and $h = h_{N,n}$ is a bandwidth which tends to zero
as $N$ and $n$ tend to infinity.

\medskip

\noi \textbf{Assumption A$_7$} \ \ $K : [-1,1] \to \R$ is a continuous function of bounded variation
that satisfies $\int_{-1}^1 K(x) \d x = 1$. Set $\| K\|^2_2 := \int_{-1}^1 K (x)^2 \d x$ and $\mu_2 (K) := \int_{-1}^1 x^2 K (x) \d x$ and $K(x) := 0$, $x \in \R \setminus [-1,1]$.

\medskip

We consider two cases separately.

\smallskip

\noi \textbf{Case (i)} $\P( b_1 = 0 ) = 1$, meaning that the coefficient $b_i=0$ for the common shock in \eqref{def:innov} is zero and that the individual processes $\{ X_i (t) \}$, $i=1,2, \ldots,$ are independent and satisfy
\begin{equation*}
X_i(t) \ = \ a_i X_i (t-1) + c_i \xi_i(t), \quad t \in \Z.
\end{equation*}

\smallskip

\noi \textbf{Case (ii)} $\P( b_1 \not = 0 ) > 0$, meaning that $\{ X_i (t) \}$, $i=1,2, \ldots,$ are mutually dependent processes.

\begin{prop}\label{prop:3}  Let Assumptions B and A$_7$ hold. If  $h^{1+p} n^{(p/2) \wedge (p-1)} \to \infty$, then
\begin{eqnarray}
\label{bias_g}
\E  \widehat{g}_{N,n} (x)  \ \to \ g (x)
\end{eqnarray}
at every continuity point $x \in \R$ of $g$. Moreover, if
\begin{eqnarray}\label{cond:Nnh}
\begin{cases}
n^{(p/2) \wedge (p-1)} h^{1+p} \to \infty & \text{ in Case {\rm (i)}},\\
n^{(p/2) \wedge (p-1)} (h/N)^{1+p} \to  \infty & \text{ in Case {\rm (ii)}},
\end{cases}
\end{eqnarray}
then
\begin{eqnarray}\label{var_g}
N h \cov ( \widehat{g}_{N,n} (x_1), \widehat{g}_{N,n} (x_2) ) \ \to \
\begin{cases}
g(x_1) \| K \|^2_2 & \text{ if } x_1 = x_2,\\
0 & \text{ if } x_1 \neq x_2
\end{cases}
\end{eqnarray}
at any continuity points $x_1, x_2 \in \R$ of $g$. If $Nh \to \infty$ holds in addition to \eqref{cond:Nnh}, then the estimator $\widehat{g}_{N,n} (x)$ is
consistent at each continuity point $x\in \R$:
\begin{equation} \label{consist}
\E | \widehat{g}_{N,n}(x) - g(x) |^2 \ \to\ 0.
\end{equation}

\end{prop}

\noi{\sc Proof.} Throughout the proof, let $K_h (x) :=  K(x/h)$, $x\in \R$. Consider \eqref{bias_g}. Note
$\E \widehat{g}_{N,n} (x)  = h^{-1} \E K_h ( x-\widehat{a}_n )$, because $\widehat{a}_{i,n}$, $i=1, \ldots, N$, are identically distributed. Let us approximate $\widehat{g}_{N,n} (x)$ by
\begin{equation}\label{def:g_N}
\widehat g_{N} (x) \ := \ \frac{1}{Nh} \sum_{i=1}^N K_h ( x-a_i ), \quad x \in \R,
\end{equation}
which satisfies $\E  \widehat g_{N} (x)  = h^{-1} \E K_h ( x-a ) \to g(x)$ as $h \to 0$ at a continuity point $x$ of $g$, see \cite{bib:PA62}. Integration by parts and
Corollary \ref{cor:1} yield
\begin{eqnarray}\label{ineq:diffE}
h | \E \widehat{g}_{N,n} (x)  - \E  \widehat g_{N} (x)  | &=&\Big|\int_{\R} (G_n(y) - G(y)) \d K_h(x-y)\Big|  \\
&\le& V(K)  \sup_{y \in [-1,1] } |G_n (y) -  G(y) | \ =  \
O ( n^{-((p/2) \wedge (p-1))/(1+p)} ), \nn
\end{eqnarray}
uniformly in $x \in \R$,  where 
$V(K)$ denotes the total variation of $K$ and $V(K) = V(K_h)$. This  proves \eqref{bias_g}.
\smallskip

Next, let us prove \eqref{var_g}. We have
\begin{eqnarray}
Nh \cov (\widehat{g}_{N} (x_1),\widehat{g}_{N} (x_2)) \ = \
\frac{1}{h}\, \E K_h ( x_1 - a ) K_h ( x_2 - a ) \ \to \
\begin{cases}
g(x_1) \| K \|^2_2 &\text{ if } x_1 = x_2,\\
0 &\text{ if } x_1 \neq x_2,
\end{cases}
\end{eqnarray}
as $h \to 0$ at any points $x_1, x_2$ of continuity of $g$, see \cite{bib:PA62}.
Split $Nh \{ \cov (\widehat{g}_{N,n} (x_1),\widehat{g}_{N,n} (x_2)) - \cov (\widehat{g}_{N} (x_1),\widehat{g}_{N} (x_2))\}
= \sum_{i=1}^3 Q_{i}(x_1,x_2)$, where
\begin{eqnarray*}
Q_1(x_1,x_2)&:=&h^{-1}\{ \E K_h ( x_1 - \widehat{a}_{n} ) K_h ( x_2 - \widehat{a}_{n} ) - \E K_h ( x_1-a ) K_h ( x_2-a )\}, \\
Q_2(x_1,x_2)&:=&h^{-1} \{ \E K_h (  x_1 -\widehat a_n  ) \E K_h ( x_2 - \widehat{a}_{n} ) -
\E K_h (  x_1 - a) \E K_h ( x_2 - a )\}, \\
Q_3(x_1,x_2)&:=&(N-1) h^{-1} \cov ( K_h ( x_1-\widehat{a}_{1,n} ) , K_h ( x_2-\widehat{a}_{2,n} ) ).
\end{eqnarray*}
Note $Q_3(x_1,x_2) = 0$ in Case (i). Similarly to \eqref{ineq:diffE},
$$
|Q_1(x_1,x_2)|
= h^{-1} \Big|\int_{\R} (G_n(y) - G(y)) \d K_h(x_1-y) K_h (x_2 -y) \Big| \le  C h^{-1} n^{-((p/2) \wedge (p-1))/(1+p)} \to 0 $$
since $ V(K_h (x_1 - \cdot) K_h(x_2 - \cdot)) \le C$
and $|Q_2(x_1,x_2)|  \le C h^{-1} n^{-((p/2) \wedge (p-1) ) /(1+p)} \to 0 $
uniformly in $x_1, x_2$. Finally,
\begin{eqnarray*}
&&\hskip-1cm|Q_3(x_1,x_2)|\\
&=& \frac{N-1}{h}
\Big|\int_{\R} \int_{\R} (\P(\widehat a_{1,n} \le y_1, \widehat a_{2,n} \le y_2) - \P(\widehat a_{1,n} \le y_1) \P(\widehat a_{2,n} \le y_2))
\d K_h(x_1-y_1)  \d K_h(x_2-y_2)   \Big|  \\
&\le& \frac{C N}{h} \sup_{y_1, y_2 \in [-1,1] } |\P(\widehat a_{1,n} \le y_1, \widehat a_{2,n} \le y_2) -
\P(\widehat a_{1,n} \le y_1) \P(\widehat a_{2,n} \le y_2)| \\
&=&O ( N h^{-1} n^{-((p/2) \wedge (p-1) )/(1+p)} ) \ = \ o(1),
\end{eqnarray*}
proving \eqref{var_g} and the proposition. \hfill $\Box$

\begin{rem} {\rm It follows from the proof of the above proposition that in the case
of a (uniformly) continuous density $g(x), x \in [-1,1]$, relations  \eqref{bias_g}, \eqref{consist}
 and the first relation in \eqref{var_g} hold uniformly in $x \in \R$, implying the convergence
 of the mean integrated squared error:
$$
\int_{-\infty}^\infty \E | \widehat{g}_{N,n}(x) - g(x) |^2 \, \d x \ \to \ 0.
$$
}
\end{rem}

\begin{prop}\label{prop:4} {\rm (Asymptotic normality)}
Let Assumptions B and A$_7$ hold and assume
$Nh \to \infty$ in addition to \eqref{cond:Nnh}.
Moreover, let $K$ be a Lipschitz function in Case {\rm (ii)}. Then
\begin{equation}\label{gclt}
 \frac{ \widehat{g}_{N,n}(x) - \E \widehat{g}_{N,n} (x)} {\sqrt{ \var (\widehat{g}_{N,n} (x)) } } \ \to_{\rm d} \ {\mathcal N}(0,1)
\end{equation}
at every continuity point $x \in (-1,1)$ of $g$.
\end{prop}

\noi{\sc Proof.} First,
consider Case (i). Since  $\widehat{g}_{N,n} (x) = (Nh)^{-1} \sum_{i=1}^N V_{i,N} $ is a (normalized) sum of i.i.d.\ r.v.s
$V_{i,N} := K_h (x-\widehat{a}_{i,n})$ with common distribution $V_N := V_{1,N}$, it suffices to verify Lyapunov's condition
\begin{eqnarray}\label{cond:lyapunov}
\frac{\E | V_N - \E V_N |^{2+\delta}}{ N^{\delta / 2} \, \{ \var (V_N) \}^{(2+\delta)/2}} \ \to \ 0,
\end{eqnarray}
for some $\delta > 0$. This follows by the same arguments as in \cite{bib:PA62}. Analogously to Proposition \ref{prop:3}, we have
$\E | V_N |^{2+\delta}  = \E | K_h ( x-\widehat{a}_n ) |^{2+\delta}  \sim h g(x) \int_{-1}^1 | K(y) |^{2+\delta} \d y = O(h)$ while
$\var (V_N) =  Nh^2 \var (\widehat g_{N,n}(x)) \sim h g(x) \|K\|^2_2 $ according to \eqref{var_g}. Hence the l.h.s. of
\eqref{cond:lyapunov} is $O( (Nh)^{-\delta/2}) = o(1)$, proving  \eqref{gclt} in Case (i).

Let us turn to Case (ii). It suffices to prove that
$\sqrt{Nh} (\widehat{g}_{N,n} (x) - \widehat g_N (x) )  \to_{\rm p} 0$, for $\widehat g_N (x)$ given in \eqref{def:g_N}.
By $|K(x) - K(y)| \le L_K | x - y |$, $x,y \in \R$, for $\epsilon >0$
\begin{eqnarray*}
\P \Big( \sqrt{Nh} | \widehat{g}_{N,n} (x) - \widehat g_N (x) | > \epsilon \Big) &\le& \P \Big(
\frac{L_K}{\sqrt{Nh}}  \sum_{i =1}^N \frac{| \widehat a_{i,n} - a_i |}{h} > \epsilon \Big)\\
&\le& N \P \Big(  | \widehat a_n - a | >
\sqrt{Nh} \Big( \frac{h}{N} \Big) \frac{\epsilon}{L_K} \Big)\\
&\le& C \Big( h (Nh)^{-p/2} \Big( \frac{N}{h} \Big)^{1+p} n^{- ((p/2) \wedge (p-1))} + \frac{N}{n} \Big) \ = \ o(1)
\end{eqnarray*}
from Proposition~\ref{prop:1} and \eqref{cond:Nnh} with $Nh \to \infty$.
\hfill $\Box$

\begin{cor}\label{cor:4} Let assumptions of Proposition \ref{prop:4} hold with $h \sim c N^{-1/5}$ for some $c>0$, i.e.,
\begin{eqnarray*}
N = \begin{cases}o (n^{ \frac{5}{3} \frac{1}{1+p} ( \frac{p}{2} \wedge (p-1) )
} ) \quad &\text{in Case {\rm (i)}},\\
o (
n^{ \frac{5}{6} \frac{1}{1+p}  ( \frac{p}{2} \wedge (p-1) ) }
) \quad &\text{in Case {\rm (ii)}.}
\end{cases}
\end{eqnarray*}
Moreover, let $g \in C^2[-1,1]$ and  $\int_{-1}^1 y K(y) \d y = 0$. Then
$$
N^{2/5} ( \widehat{g}_{N,n} (x) - g(x) ) \ \to_{\rm d} {\mathcal N} ( \mu(x), \sigma^2(x) ),
$$
where $\mu(x) := (c^2/2) g''(x) \mu_2 (K)$ and $\sigma^2(x) := (1/c) g(x) \| K \|^2_2$.
\end{cor}

\noi{\sc Proof.} This follows from Proposition \ref{prop:4}, by noting that
$\E \widehat g_{N} (x)  - g(x) \sim h^2 g''(x) \mu_2 (K) / 2$ as $h \to 0$ and $\E  \widehat{g}_{N,n} (x)  - \E  \widehat g_N (x)  = O (h^{-1} n^{-((p/2) \wedge (p-1))/(1+p)})$ by \eqref{ineq:diffE}.
\hfill $\Box$

\section{Simulations
}\label{s:sim}

In this section we compare our nonparametric goodness-of-fit test in \eqref{gof} for  testing the null hypothesis $G = G_0$
with its parametric analogue studied in \cite{bib:BSG10}. In accordance with the last paper, we
assume $\{ X_i (t) \}$ in \eqref{def:multple_AR1} to be independent AR(1) processes with standard normal i.i.d.\ innovations $\{\zeta_i (t)
\}$,  ${\zeta(0) \sim {\mathcal N}(0,1)}$ and
the random autoregressive coefficient $a_i \in (0,1)$ having a Beta-type density  $g(x)$
with unknown parameters $\theta := (\alpha, \beta)^\top$:
\begin{eqnarray}\label{def:pdf}
g(x) \ = \ \frac{2}{{\rm B}(\alpha, \beta)} x^{2 \alpha - 1} (1 - x^2)^{\beta - 1}, \quad x \in (0,1), \quad \alpha > 1, \beta > 1.
\end{eqnarray}
Note that $\beta\in (1,2)$ implies the long memory property in
$\{X_i(t)\}$. Beran et al.\
\cite{bib:BSG10} discuss a maximum likelihood estimate $\widehat{\theta}_{N,n,\kappa} = (\widehat{\alpha}, \widehat{\beta})^\top$ of $\theta= (\alpha, \beta)^\top$ when each unobservable coefficient $a_i$ is replaced by its estimate $\widehat{a}_{i, n, \kappa} := \min \{ \max \{ \widehat{a}_{i,n}, \kappa \}, 1-\kappa \}$
with $\widehat{a}_{i,n}$ given in \eqref{hatan} and $0< \kappa = \kappa(N,n) \to 0 $ is a truncation parameter.
Under certain conditions on $N, n \to \infty$ and $\kappa \to 0$, Beran et al.~\cite[Theorem 2]{bib:BSG10} showed that
\begin{equation}\label{T2}
N^{1/2} ( \widehat{\theta}_{N,n,\kappa} - \theta_0) \ \to_{\rm d} \ {\mathcal N}(0, A^{-1} (\theta_0)),
\end{equation}
where $\theta_0$ is the true parameter vector,
\begin{eqnarray*}
A (\theta) :=
\begin{pmatrix}
\psi_1(\alpha) - \psi_1(\alpha+\beta) & - \psi_1(\alpha+\beta)\\
 - \psi_1(\alpha+\beta) & \psi_1(\beta) - \psi_1(\alpha+\beta)
\end{pmatrix},
\end{eqnarray*}
and $\psi_1(x) := \d^2 \ln \Gamma (x)/\d x^2$ is the Trigamma function. Based on \eqref{gof} and \eqref{T2}, we consider testing
both ways (nonparametrically  and parametrically) the hypothesis that
the unobserved autoregressive coefficients $a_1, \ldots, a_N$ are drawn from  the reference distribution $G_0$
having density function in \eqref{def:pdf} with a specific $\theta_0$, i.e., the null
 $G = G_0$ vs.\ the alternative $G \not = G_0$. The respective test statistics are
\begin{eqnarray} \label{TT}
T_1 :=   N^{1/2} \sup_x | \widehat{G}_{N,n}(x) -G_0(x) | \quad \text{and}  \quad
T_2 :=  N (\widehat{\theta}_{N,n,\kappa} - \theta_0)^\top A(\theta_0) (\widehat{\theta}_{N,n,\kappa} - \theta_0).
\end{eqnarray}
Under the null hypothesis, the distributions of statistics $T_1$  and $T_2$ converge to the Kolmogorov distribution
and the chi-square distribution with 2 degrees of freedom, respectively, see \eqref{gof}, \eqref{T2}.

To compare the performance of the above testing procedures, we compute the empirical
distribution of the p-value of $T_1$ and $T_2$ under null and alternative hypotheses. The p-value of observed $T_i$
is defined as $p(T_i) = 1 - {\cal K}_i(T_i), i=1,2,  $ where ${\cal K}_i(y), i=1,2 $  denote the
limit distribution functions of \eqref{TT}. Recall that when the significance level of the test is correct,  the (asymptotic)
distribution of the p-value is uniform on $[0,1]$.  The simulation procedure to compare the performance of $T_1$ and $T_2$ is the following:
\smallskip

\noi \textbf{Step S$_0$} \ We fix the parameter under the null hypothesis $H_0 : \theta= \theta_0$ with   $\theta_0 = (2, 1.4)^\top$.

\smallskip

\noi \textbf{Step S$_1$} \ We simulate $5000$ panels with $N=250$, $n=817$ for five chosen values $\theta = (2, 1.2)^\top,  (2, 1.3)^\top, $  $
(2, 1.4)^\top, (2, 1.5)^\top, (2, 1.6)^\top$ of Beta parameters.

\smallskip

\noi \textbf{Step S$_2$} \  For each simulated panel we compute the p-value of statistics $T_1$ and $T_2$.

\smallskip

\noi \textbf{Step S$_3$} \  The empirical c.d.f.'s of computed p-values of statistics $T_1$ and $T_2$
are plotted.

\smallskip

The values of Beta parameters $\theta_0 = (2, 1.4)^\top$, $N$, $n$ were chosen in accordance with the simulation study in \cite{bib:BSG10}.

Fig.\ 1 presents the simulation results under the true hypothesis $\theta = \theta_0$ with zoom-in on small p-values. We see that both c.d.f.'s in the left graph
are approximately linear.
Somewhat surprisingly, it appears that the empirical size of $T_1$ (the nonparametric test) is better
than the size of $T_2$ (the parametric test). Particularly, for significance levels $0.05$ and $0.1 $ we provide the empirical size values in Table 1.

Fig.\ 2 gives the graphs of the empirical c.d.f.'s of p-values of $T_1$ and $T_2$
for several alternatives $\theta \ne \theta_0$. It appears that for $\beta > \beta_0 = 1.4 $
the parametric test  $T_2$ is more powerful than the nonparametric test $T_1$
but for $\beta < \beta_0$ the power differences are less significant. Table 1 illustrates the empirical power for the significance levels $0.05$, $0.1$.

\begin{figure}[h]
\centering
\includegraphics[width=15cm]{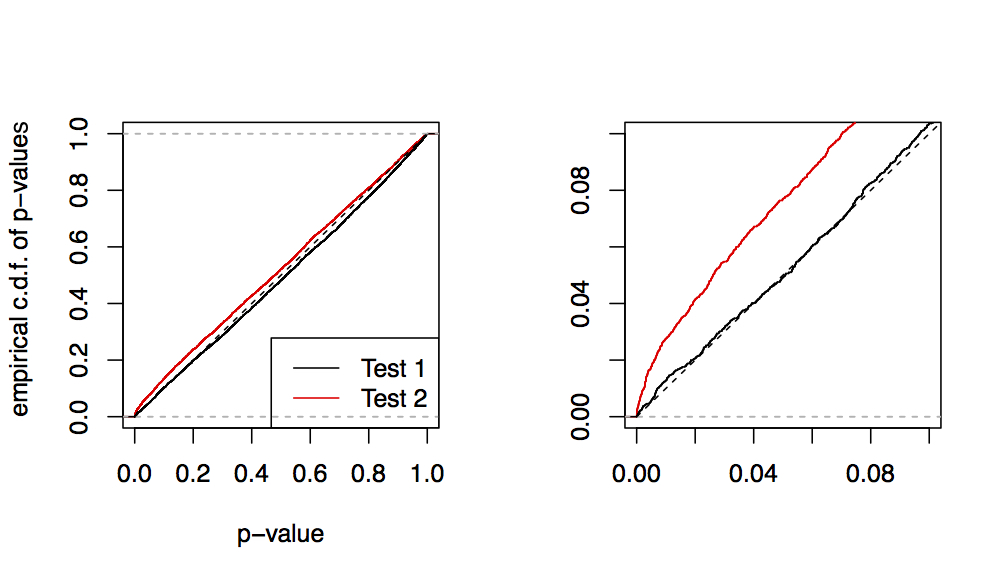}
\caption{[left] Empirical c.d.f. of p-values of $T_1$ and $T_2$
under $H_0: \theta_0  = (2, 1.4)^\top$; 5000 replications with $N=250$, $n=817$. [right] Zoom-in on the region of interest: p-values smaller than 0.1.}
\end{figure}

\begin{figure}[h]
  \centering
  \includegraphics[width=15cm]{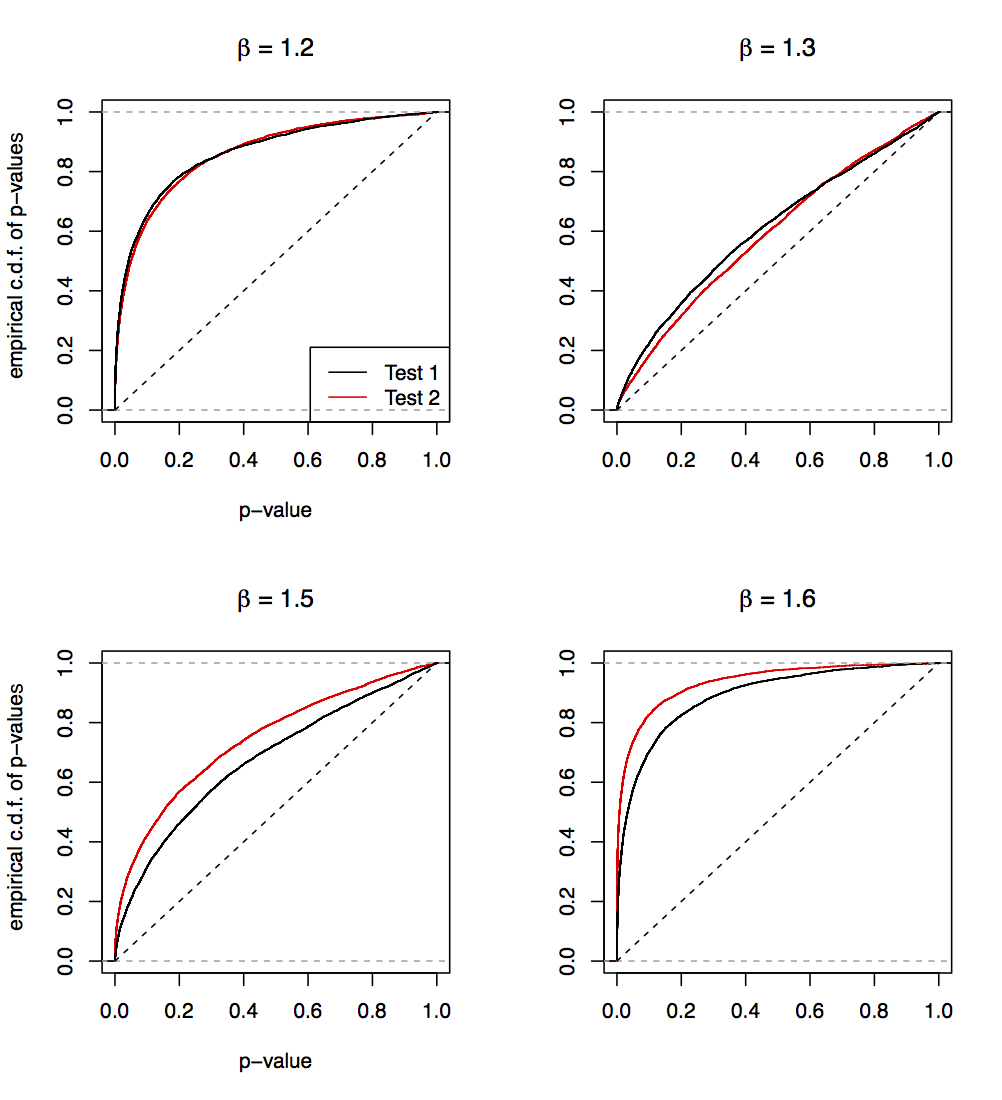}
\caption{Empirical c.d.f. of p-values of $T_1$ and $T_2$ for testing $H_0 : \theta_0 = (2,1.4)^\top$
under several alternatives of the form $\theta = (2, \beta)^\top$; 5000 replications with $N=250$, $n=817$.}
\end{figure}

\begin{table}[htdp]
\begin{center}
\begin{tabular}{c|ccccc|ccccc}
\hline
Signif. level  & \multicolumn{5}{|c|}{ 5\% } & \multicolumn{5}{|c}{ 10\% } \\
\hline
 $\beta$  &  1.2 & 1.3 &  1.4 & 1.5 & 1.6 &  1.2 & 1.3 &  1.4 & 1.5 & 1.6\\
\hline
$T_1$   & .532 & .137 & .049 & .208 & .576 & .653 & .223 & .103  & .315 & .702 \\
 $T_2$   &  .500 & .104 & .077 & .313 & .735  & .634 & .184 & .134 & .421 &.827 \\
\hline
\end{tabular}
\caption{ Numerical results  of the comparison for testing procedure $ H_0\colon \theta = (2,1.4)^\top$ at significance level 5\% and 10\% . The column for $\beta=1.4$ provides the empirical size.}\label{tab:1}
\end{center}
\end{table}

The above simulations (Fig.\ 1 and 2, Table 1) refer to the case of independent individual processes  $\{ X_i (t) \}$. There are no theoretical results for the parametric test $T_2$, when AR(1) series are dependent. Although the nonparametric test $T_1$ is valid for the latter case,  one may expect
that the presence of the common shock component in the panel data in \eqref{def:innov} has a negative effect on the test  performance for short series. To illustrate this effect, we simulate 5000 panels with AR(1) processes  $\{ X_i (t) \}$ driven by dependent shocks in \eqref{def:innov} with $b_i = b$, $c_i = (1 - b^2)^{1/2}$. As previously, we choose $\theta_0 = (2, 1.4)^\top$, $N=250$, $n=817$ and  we fix $\theta = (2, 1.4)^\top$ to evaluate the empirical size of $T_1$.
Fig.\ 3 [left] presents the graphs of the empirical c.d.f.'s of the p-values of $T_1$
for
$b = 1$, $b =0.6$ and $b=0$, the latter corresponding to independent individual processes as in Fig.\ 1.
We see that the size of the test worsens when $b$ increases, particularly when $b=1$ and  the individual processes
are all driven by the same
common noise. To overcome the last effect, the sample length $n$ of each series in the panel may be increased as in
Fig.\ 3 [right], where the choice of $n = 5500$ and $b=1$ shows a much better performance
of $T_1$ under the null hypothesis $\theta = \theta_0 = (2,1.4)^\top $ and the alternative ($\theta = (2, 1.5)^\top$  and $\theta = (2,1.6)^\top$) scenarios.

\begin{figure}[t!]
  \centering
  \includegraphics[width=15cm]{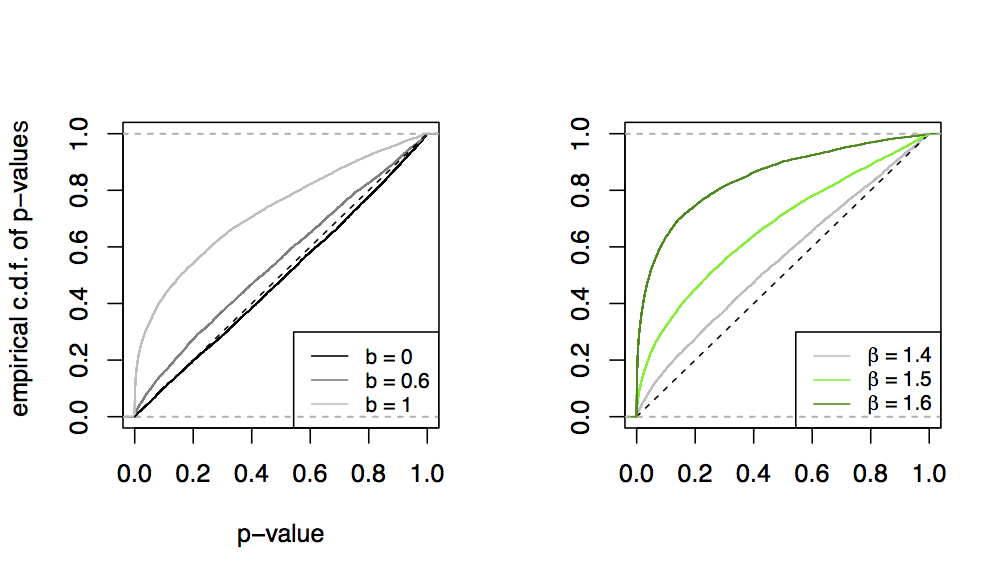}
\caption{[left] Empirical c.d.f.\ of p-values of $T_1$ under $H_0\colon \theta_0  = (2, 1.4)^\top$ for different dependence structure between AR(1) series
: $b_i = b$ and $c_i = \sqrt{1-b^2}$ and $N=250$, $ n=817$.
[right] Empirical c.d.f. of p-values of $T_1$ for testing $H_0 :  \theta_0  = (2, 1.4)^\top$. AR(1) series are driven by common innovations, i.e., $b_i = 1$, $c_i = 0$, for $\theta  = (2, \beta)^\top$; 5000 replications with $N=250$, $n = 5500$.
}\label{fig:T1_dep}
\end{figure}

\bigskip

In conclusion,
\begin{enumerate}
\item We do not observe an important loss of the power for the nonparametric KS test $T_1$ compared to the parametric approach.

\item The KS test $T_1$ does not require to choose any tuning parameter contrary to the test $T_2$.

\item One can use the KS test $T_1$ under weaker assumptions on AR(1) innovations. We only impose moment conditions. The dependence between the series is allowed by \eqref{def:innov}.
\end{enumerate}

\section*{Acknowledgements}
We thank the referees for valuable comments that led to an improved version of this paper. The first, third and fourth authors also acknowledge the support by a grant (No.~MIP-063/2013) from the Research Council of Lithuania.

\clearpage

\bibliographystyle{plain}

\input bibnames.sty

\section{Appendix: some proofs and auxiliary lemmas
}

We use the following martingale moment inequality.

\begin{lemma} \label{Yp} Let $p>1 $ and
$\{\xi_j, j \ge 1\}$ be  a martingale difference sequence:
$\E [\xi_j |\xi_1, \ldots, \xi_{j-1}] = 0$, $j=2,3, \ldots $ with $\E |\xi_j|^p < \infty $.
Then there exists a constant $C_p<\infty$ depending only on $p$ and such that
\begin{equation}\label{ineq:pmom}
\E \Big|\sum_{j=1}^\infty \xi_j\Big|^p \ \le \ C_p \begin{cases}
\sum_{j=1}^\infty \E |\xi_j|^p, &1< p \le 2, \\
\big(\sum_{j=1}^\infty (\E |\xi_j|^p)^{2/p}\big)^{p/2}, &p > 2.
\end{cases}
\end{equation}
\end{lemma}

For $1 < p \le  2 $, inequality \eqref{ineq:pmom} is known as von Bahr and Ess\'een inequality, see \cite{bib:BE65},
and for $p> 2$,  it is
a consequence of the Burkholder and Rosenthal inequality (\cite{bib:BU73, bib:RO70}, see also  \cite[Lemma~2.5.2]{bib:book2012}).

\medskip

\noi{\sc Proof of Proposition \ref{prop:1}. } Since $\widehat a_n$ in \eqref{han} is invariant w.r.t.\ a scale factor
of innovations $\{ \zeta (t) \}$, w.l.g.\ we can assume $b^2 + c^2 =1 $ and
$\E \zeta^2 (0) = 1, $  $\E | \zeta (0) |^{2p} < \infty$.
Then $\widehat{a}_n - a  = \sum_{i=1}^3 \delta_{ni}$, where
\begin{eqnarray*}\label{def:deltan1}
\delta_{n1} &:=& -\frac{a X^2(n)}{\sum_{t=1}^n X^2(t)-  n (\bar X_n)^2}, \quad \delta_{n2} \ := \  \frac{\sum_{t=1}^{n-1}X(t) \zeta(t+1)
}{\sum_{t=1}^{n} X^2(t)-  n (\bar X_n)^2},\\
\delta_{n3} &:=&  \frac{\bar X_n (X(1) + X(n)) - (\bar X_n)^2 (1 + n(1-a))}
{\sum_{t=1}^{n} X^2(t)- n (\bar X_n)^2}.
\end{eqnarray*}
The statement of the proposition follows from
\begin{equation}\label{delta12}
\P( | \delta_{ni} | > \gamma) \le C(n^{-1} + n^{-(p/2) \wedge (p-1)} \gamma^{-p}) \qquad (0< \gamma < 1, \, i=1,2,3).
\end{equation}
To show \eqref{delta12} for $i=1$, note that $\delta_{n1} = L_n/(n + D_n)$,
where $L_n := - a (1-a^2) X^2(n)$ and $D_n = D_{n1} - D_{n2}$, $D_{n1} := \sum_{t=1}^n ((1-a^2) X^2(t) -1)$, $D_{n2} := n (1-a^2) (\bar X_n)^2$.
We have $\P( | \delta_{n1} | > \gamma) \le  \P( |D_{n}| > n/2) + \P( |L_{n}| > n \gamma/2) $.
Thus,   \eqref{delta12} for $i=1$
follows from
\begin{eqnarray} \label{S12}
\E  |D_{n1}|^{p \wedge 2} &\le&Cn, \qquad \E |D_{n2}| \ \le \ C \qquad \text{and}  \qquad
\E |L_{n}|^p \ \le\  C.
\end{eqnarray}

Consider the first relation in \eqref{S12}.
Clearly, it suffices to prove it for $1< p \le 2 $ only.
We have $D_{n1} = 2 D'_{n1}+ D''_{n1}$,
where
\begin{eqnarray*}
D'_{n1}&:=&  (1-a^2) \sum_{s_2 < s_1 \le n}  \sum_{t=1 \vee s_1}^n a^{2(t-s_1)} a^{s_1-s_2} \zeta (s_1) \zeta(s_2), \\
D''_{n1}&:=&(1- a^2)\sum_{s \le n} \sum_{t=1\vee s}^n a^{2(t-s)} (\zeta^2(s)-1).
\end{eqnarray*}
We will use the following elementary inequality: for any  $-1 \le a \le 1, n \ge 1, s \le n $
\begin{eqnarray}
\alpha_n (s)\ :=\ (1-a^2) \sum_{t=1\vee s}^n a^{2(t-s)}
&=&\begin{cases}
 a^{2(1-s)} (1 - a^{2n}), &s \le 0, \\
 1- a^{2(n + 1 - s)}, &1 \le s \le n
         \end{cases} \nn \\
&\le&C \begin{cases}
a^{-2s} \min(1, 2n (1 - |a|)), &s \le 0, \\
1, &1 \le s \le n.
\end{cases}     \label{alphain}
\end{eqnarray}
Using the independence of $\{\zeta(s)\}$ and $a $
and inequality \eqref{ineq:pmom} (twice) for $1< p\le 2 $ we obtain
\begin{eqnarray*}
\E |D'_{n1} |^p
&=&\E \big|\sum_{s_1 \le n} \alpha_n(s_1) \zeta(s_1) \sum_{s_2<s_1} a^{s_1-s_2} \zeta(s_2)\big|^p \\
&\le&C\E \sum_{s_1 \le n}  \big| \alpha_n (s_1) \zeta (s_1) \sum_{s_2 < s_1} a^{s_1-s_2} \zeta (s_2) \big|^p\\
&\le&C\E \sum_{s_1 \le n}  |\alpha_n (s_1)|^p  \sum_{s_2 < s_1} |a|^{p(s_1-s_2)}  \\
&\le&C \E (1- |a|)^{-1} \sum_{s \le n} |\alpha_n(s)|^p  \ \le \ Cn
\end{eqnarray*}
since  $ \E (1- |a|)^{-1} < \infty $ (see \eqref{cond:Ea}) and
$\sum_{s \le n} |\alpha_n(s)|^p \le C n$ follows from \eqref{alphain}.
Similarly, since $\{ \zeta^2 (s) - 1, \, s \le n \}$ form a martingale difference sequence,
\begin{eqnarray*}
 \E |D''_{n1} |^p
 &\le&C\E  \sum_{s\le n} |\alpha_n(s)|^p  \ \le \  Cn,
\end{eqnarray*}
proving the first inequality \eqref{S12}.
The second inequality in \eqref{S12} follows by noting that $n \bar X_n = \sum_{s \le n} ( \sum_{t=1 \vee s}^n a^{t-s} ) \zeta (s)$
and
$$
(1-a^2) \E [ (n \bar X_n)^2 | a] = a^2 \Big( \frac{1-a^n}{1-a} \Big)^2 + (1-a^2) \sum_{s=1}^n \Big( \frac{1-a^s}{1-a} \Big)^2 \le \frac{Cn}{1-a}.
$$
Consider the last inequality in \eqref{S12}. We have $| L_{n} | \le |2 L'_{n}+ L''_{n} + 1|$,
where
\begin{eqnarray*}
L'_{n}&:=&(1-a^2) \sum_{s_2 < s_1 \le n}  a^{2(n-s_1)} a^{s_1-s_2} \zeta (s_1) \zeta(s_2), \qquad
L''_{n}\ :=\ (1- a^2)\sum_{s \le n} a^{2(n-s)} (\zeta^2(s)-1).
\end{eqnarray*}
We use Lemma \ref{Yp}, as above.
Let  $1 \le p \le 2 $. Then
$ \E |L''_{n}|^p
\le C \E \sum_{s \le n} \{ (1-a^2) a^{-2(n-s)} \}^p \le  C $ and
$\E |L'_{n}|^p
\le C \E \sum_{s_2 < s_1 \le n} \{(1-a^2) |a|^{2(n-s_1)} |a|^{s_1-s_2}\}^p \le  C \E (1-|a|)^{p-2} \le C$.
Next, let $p \ge 2 $. Then $ \E |L''_{n}|^p
\le C \E \{ \sum_{s \le n} |(1-a^2) a^{2(n-s)}|^2 \}^{p/2} \le C $ and
$ \E |L'_{n}|^p
\le C \E (1-a^2)^p \{\sum_{s_2 < s_1 \le n} a^{4(n-s_1)} a^{2(s_1-s_2)} \}^{p/2} \le C$,
proving \eqref{S12} and hence \eqref{delta12} for $i=1$.

Consider \eqref{delta12} for $i=2$.  We have $\delta_{n2} = R_{n}/(n + D_{n})$,
where $R_{n} := (1-a^2)
 \sum_{t=1}^{n-1} X(t) \zeta (t+1)$ and  $D_{n} $ 
is the same as in
\eqref{S12}. Then
$\P(|\delta_{n2}| > \gamma) \le  \P(|R_{n}| > n\gamma/2) + \P(|D_{n}| > n/2) $,
where
\begin{eqnarray}\label{S1}
\P(|D_{n}| > n/2) &\le& (n/4)^{-(p\wedge 2)}  \E |D_{n1}|^{p \wedge 2} + (n/4)^{-1} \E |D_{n2}|\nn \\
&\le&
C\begin{cases} n^{-(p-1)}, &1< p \le 2, \\
n^{-1}, &p > 2,
\end{cases}
\end{eqnarray}
according to  \eqref{S12}. Therefore
\eqref{delta12} for $i=2$ follows
from
\begin{equation} \label{R}
\E  |R_{n}|^p \le C\begin{cases}
n, &1< p \le 2, \\
n^{p/2}, &p>2.
\end{cases}
\end{equation}
Since
$R_n = (1 - a^2)\sum_{s  \le n-1} \zeta (s) \sum_{t=1\vee s}^{n-1} a^{t-s} \zeta (t+1) $ is a sum of martingale differences,
by inequality \eqref{ineq:pmom} with $1< p \le 2 $
we obtain
\begin{eqnarray*}
\E |R_n|^p
&\le&C\E \sum_{s \le n-1} \big| (1-a^2) \zeta (s) \sum_{t=1\vee s}^{n-1} a^{t-s} \zeta (t+1) \big|^p \\
&\le&C\E |1-a^2|^p  \sum_{s \le n-1} \sum_{t=1\vee s}^{n-1}|a|^{p(t-s)} \\
&\le&C \E |1-a^2|^p \big( \sum_{s\le 0} |a|^{-ps} \sum_{t=1}^{n-1} |a|^{pt}
+ \sum_{s=1}^{n-1} \sum_{t=s}^{n-1} |a|^{p(t-s)} \big) \\
&\le&C \E |1-a^2|^p \big\{ (1-|a|^p)^{-2} + n (1- |a|^p)^{-1} \big\}  \ \le  Cn,
\end{eqnarray*}
proving \eqref{R} for $p\le 2$. Similarly, using \eqref{ineq:pmom} with $p>2 $ we get
\begin{eqnarray*}
\E |R_n|^p
&=& \E \Big[  |1-a^2|^p \E \big[
| \sum_{s \le n-1} \zeta (s) \sum_{t=1\vee s}^{n-1} a^{t-s} \zeta (t+1) |^p \big| a \big] \Big]\\
&\le&C\E \Big[ |1-a^2|^p  \big\{ \sum_{s \le n-1} \big(\E \big[ | \zeta(s) \sum_{t=1\vee s}^{n-1}a^{t-s} \zeta(t+1) |^p \big| a\big]\big)^{2/p} \big\}^{p/2} \Big] \\
&\le&C \E |1-a^2|^p \big\{ \sum_{s \le n-1} \sum_{t=1\vee s}^{n-1} a^{2(t-s)}  \big\}^{p/2} \\
&\le&C \E |1-a^2|^p \big\{ \sum_{s\le 0} a^{-2s} \sum_{t=1}^{n-1} a^{2t}
+ \sum_{s=1}^{n-1} \sum_{t=s}^{n-1} a^{2(t-s)} \big\}^{p/2} \\
&\le&C \E |1-a^2|^p \big\{ (1-a^2)^{-2} + n (1- a^2)^{-1} \big\}^{p/2} \
\le \ Cn^{p/2},
\end{eqnarray*}
proving \eqref{R} and \eqref{delta12} for $i=2$.

It remains to prove \eqref{delta12} for $i=3$.  Similarly as above,
$\P(|\delta_{n3}| > \gamma) \le  \P(|Q_{n}| > n\gamma/2) + \P(|D_{n}| > n/2) $, where
$Q_n := (1-a^2) \{\bar X_n (X (1) + X(n)) - (\bar X_n)^2 (1 + n(1-a)) \}$ and $D_{n}$
is evaluated in \eqref{S1}. Thus, \eqref{delta12} for $i=3$ follows from \eqref{S1} and
\begin{equation}\label{Qn}
\E |Q_n|^p  \ \le \ C \{ \E |(1-a^2) X^2(n)|^p + \E |(1-a^2)(\bar X_n)^2|^p + n^p \E |(1-a)(1-a^2)(\bar X_n)^2|^p \}\ \le \ C.
\end{equation}
Since $n \bar X_n = \sum_{s \le n} (\sum_{t=1 \vee s}^{n} a^{t-s}) \zeta (s)$, an application of
the second inequality of \eqref{ineq:pmom} yields
\begin{eqnarray*}
\E [ |n \bar X_n|^{2p} | a ] \le C \Big(  \frac{(1-a^n)^{2}}{(1-a^2)(1-a)^{2}} + \sum_{s=1}^{n} \Big( \frac{1-a^s}{1-a} \Big)^2 \Big)^p.
\end{eqnarray*}
Using $1 - a^n \le 1 \wedge (n (1-a))$ we obtain
$\E | (1-a)(1-a^2)(\bar X_n)^2 |^p \le C n^{-p}$ and $\E | (1-a^2) (\bar X_n)^2 |^p \le C \E (1-a)^{-1} n^{-1}$. Finally, $\E |(1-a^2) X^2(n)|^p \le C$ follows by the same arguments as $\E |L_n|^p \le C$ (see \eqref{S12}).
This proves \eqref{Qn}, thereby completing
the proof of \eqref{delta12} and of the proposition, too. \hfill $\Box$

\medskip

Let $a, a_1, \ldots, a_N$ be i.i.d.\ r.v.s
with d.f.\  $G(x) = \P (a \le x)$ supported by $[-1,1]$.
Define $\widehat G_N (x) := N^{-1} \sum_{i=1}^N \1(a_i \le x)$, $U_N(x) := N^{1/2} (\widehat G_N(x) - G(x))$, $x \in \R$, and
$\omega_N (\delta)$ (= the modulus of continuity of $U_N$) by
$$
\omega_N (\delta) \ := \ \sup_{0 \le y-x \le \delta} | U_N (y) - U_N (x) |, \quad \delta > 0.
$$

\begin{lemma}\label{lemma:2}
Assume that $G$ satisfies Assumption A$_6$. Then for all $\epsilon > 0$,
$$
\epsilon^4 \P ( \omega_N(\delta) > 6 \epsilon ) \ \le \ (3 + 3 C) L_G \delta^\varrho + N^{-1},
$$
where $C$ is a constant independent of $\epsilon, \, \delta, \, N$.
\end{lemma}

\noi{\sc Proof.}
As in \cite[p.~106,~(13.17)]{bib:BILL68} we have that
\begin{eqnarray*}
\E | U_N (y) - U_N (x) |^2 | U_N (z) - U_N (y) |^2
&\le& 3 \P(a \in (x,y]) \P(a \in (y,z]),\\
\E | U_N (y) - U_N (x) |^4 &\le& 3 \P(a \in (x,y])^2 + N^{-1} \P(a \in (x,y])
\end{eqnarray*}
for $-1 \le x \le y \le z \le 1$, where the second inequality treats the 4th central moment of a binomial variable.
Now fix $\delta > 0$ and split $[-1, 1] = \cup_i \Delta_{i}$, where
$\Delta_{i} = [-1+ i\delta, -1 + (i+1)\delta]$, $i=0,1, \ldots, \lfloor 2/\delta\rfloor-1$, $\Delta_{\lfloor 2/\delta\rfloor} = [-1 + \lfloor 2/\delta\rfloor\delta, 1]$.
According to \cite[p.~49, Lemma~1]{bib:SW86}, for all $\epsilon > 0$,
$$
\epsilon^4 \P ( \omega_N(\delta) > 6 \epsilon ) \ \le \ (3 + 3 C) \max_i \P (a \in \Delta_i) +  N^{-1},
$$
where $C$ is a constant independent of $\epsilon, \, \delta, \, N$. Lemma follows from Assumption A$_6$ on the d.f.\ $G$ of the r.v.\ $a$.
\hfill $\Box$

\smallskip

Note that if we take $\delta = \delta_N = o(1)$,
we then get $\P ( \omega_N (\delta) > \epsilon ) \to 0$ as $N \to \infty$.

\bigskip

\begin{lemma}\label{lemma:5}
Let $\widehat a_{1,n}, \widehat{a}_{2,n}$ be given in \eqref{hatan} under Assumptions A$_1$--A$_6$ with $\varrho = 1$. Then for all $\gamma \in (0,1)$ and $n \ge 1$, it holds
\begin{equation*}
\sup_{x,y \in [-1,1]} |\P(\widehat a_{1,n} \le x, \, \widehat a_{2,n} \le y)
- \P(\widehat a_{1,n} \le x) \P(\widehat a_{2,n} \le y)| \ = \ O (n^{-((p/2) \wedge (p-1) )/(1+p)}).
\end{equation*}
\end{lemma}

\noi{\sc Proof.} Define $\delta_{i,n} := \widehat a_{i,n} - a_i$, $i =1,2$. For $\gamma \in (0,1)$, we have
\begin{eqnarray*}
\P(|\delta_{1,n}|> \gamma \, \text{or} \, |\delta_{2,n}|> \gamma ) \ \le \ \P(|\delta_{1,n}|> \gamma ) +  \P(|\delta_{2,n}|> \gamma) \
\le \  C ( n^{-((p/2) \wedge (p-1))} \gamma^{-p} + n^{-1} )
\end{eqnarray*}
by Proposition \ref{prop:1}.
Consider now
\begin{eqnarray*}
\P(\widehat a_{1,n} \le x, \, \widehat a_{2,n} \le y)
&=&\P(a_1 + \delta_{1,n} \le x, \,  a_2 + \delta_{2,n} \le y) \\
&\le&\P(a_1+ \delta_{1,n} \le x, \, a_2 + \delta_{2,n} \le y, \, |\delta_{1,n}| \le \gamma, \, |\delta_{2,n}| \le \gamma )
+ \P(|\delta_{1,n}|> \gamma \, \text{or} \, |\delta_{2,n}|> \gamma ).
\end{eqnarray*}
Then
\begin{eqnarray*}
\P(a_1 + \delta_{1,n} \le x, \, a_2 + \delta_{2,n} \le y, \, |\delta_{1,n}| \le \gamma, \, |\delta_{2,n}|\le \gamma )
&\le&\P(a_1 \le x+ \gamma,  \, a_2  \le y +\gamma, \, |\delta_{1,n}| \le \gamma, \, |\delta_{2,n}| \le \gamma )\\
&\le&
G(x+\gamma)G(y+\gamma)
\end{eqnarray*}
and
\begin{eqnarray*}
\P(a_1 + \delta_{1,n} \le x, \,  a_2 + \delta_{2,n} \le y, \, |\delta_{1,n}| \le \gamma, \, |\delta_{2,n}|\le \gamma )
&\ge&\P(a_1 \le x-\gamma, \, a_2  \le y -\gamma, |\delta_{1,n}| \le \gamma, \, |\delta_{2,n}|\le \gamma )\\
&\ge&G(x-\gamma)G(y-\gamma) - \P(|\delta_{1,n}|> \gamma \, \text{or} \, |\delta_{2,n}|> \gamma ).
\end{eqnarray*}
From \eqref{cond:GHoe} we obtain
\begin{equation*}
|G(x \pm \gamma) G(y \pm \gamma) - G(x)G(y)| \ = \ |(G(x) + O(\gamma))(G(y) + O(\gamma)) - G(x) G(y)|
\ \le \ C\gamma.
\end{equation*}
Hence,
\begin{eqnarray}\label{var1}
|\P(a_{1} \le x, \, a_{2} \le y)- G(x)G(y)| \ \le \ C ( \gamma + n^{-1} + n^{-((p/2) \wedge (p-1))} \gamma^{-p} ).
\end{eqnarray}
In a similar way,
\begin{eqnarray}\label{var2}
|\P(a_{1} \le x) \P(a_{2} \le y)- G(x)G(y)| \ \le \ C (\gamma + n^{-1} + n^{-((p/2) \wedge (p-1))} \gamma^{-p} ).
\end{eqnarray}
By \eqref{var1}, \eqref{var2}, the proof of the lemma is complete with $\gamma = \gamma_n = o(1)$, which satisfies $\gamma_n \sim n^{-(p/2) \wedge (p-1)} \gamma^{-p}_n$.
\hfill $\Box$

\end{document}